# Reward-penalty Mechanism for Reverse Supply Chain Network with Asymmetric Information and Carbon Emission Constraints

Xiao-qing Zhang[1]    Xi-gang Yuan[2]

(1. The School of Statistics Southwestern University of Finance and Economics        Chengdu 611130, P.R. China

(2. The School of Statistics Southwestern University of Finance and Economics        Chengdu 611130, P.R. China)

**Abstract:** In this paper, we discuss the government's reward and penalty mechanism in the presence of asymmetric information and carbon emission constraint when downstream retailers compete in a reverse supply chain network. Considering five game models which are different in terms of the coordination structure of the reverse supply chain network and power structure of the reward-penalty mechanism: (1) the reverse supply chain network centralized decision-making model; (2) the reverse supply chain network centralized decision-making model with carbon emission constraint; (3) the retailers' competition reverse supply chain network decentralized decision-making model; (4) the retailers' competition reverse supply chain network decentralized decision-making model with carbon emission constraint; and (5) the retailers' competition reverse supply chain network decentralized decision-making model with carbon emission constraint and the government's reward-penalty mechanism. Building the participation - incentive contract under each model use the principal-agent theory, and solving the model use the Lagrange multiplier method. We can get the following conclusion: (1) when the government implements the reward-penalty mechanism for carbon emission and recycling simultaneously, the recycling rate as well as the buy-back price offered by the manufacturer are higher than those when the government conducts reward-penalty mechanism exclusively for carbon emission; (2) when the government implements carbon emission constraint，both retailers' selling prices of the new product are higher than those when no carbon emission constraint is forced; (3), there is no certain relationship between the two retailers' selling prices of the new product when the government implements the reward-penalty mechanism only for carbon emission and when it implements the mechanism for carbon emission as well as recycling; (4) Regardless of the retailer's fixed cost is high or low, when the government implement reward-penalty mechanism for carbon emission, the retailer's buy-back price is affected by the carbon emission in the manufacturer; (5) When it satisfy certain some conditions, the buy-back price of the retailers' competition reverse supply chain decentralized decision-making model with carbon emission constraint and the government's reward-penalty mechanism can larger than that of the retailers' competition reverse supply chain decentralized decision-making model with carbon emission constraint.

## 1 Introduction

With the development of the economy and society, more and more people pay much attention to the environmental protection and ecological balance problem. There are more than 700 kinds of chemical composition in the waste electrical and electronic products (WEEE), most of them are harmful for people, on the other hand, the WEEE can bring much more economic value. Thus, the government makes much more policies and regulations so that it can guide the manufacturer or the retailer to recycle and reuse the WEEE. Many measures reflect the ideas of the extended production responsibility, such as the recycling of waste product, the standardization of recycling and manufacturing, and the batch of production. In 2003, the European Union publishes the regulation about recycling the WEEE products so that it effectively promotes the recycling of the WEEE products, and controls the environment problem.

In China, some regulations about the WEEE products recycling begin to implement. In addition, in order to implementing the extended production responsibility, manufacturer should consider the environmental problem in the processing of design, production and sale. The government policy plays an important role in the WEEE products collection. Through extend production responsibility law, the government incentivizes the manufacturer or the retailer to produce the environmentally friendly product and recycles the WEEE products at the end of the product life. In Europe, producers take responsibility for recycling the used vehicles, packaging and other WEEE products. Until now, many US states have passed law mandating state wide e-waste recycling.

From above analysis, we find that some regulations in European Union stress penalty, on the other hand, the regulation in china focuses on subsidy and penalty. We have not found the mechanism that puts together reward with penalty. In this paper, the reward and penalty mechanism is considered. With asymmetric information and retail-level competition in the presence of carbon emission constraints, we mainly identify the best design parameters for the reward and penalty mechanism so that it promotes the retailer to collecting more WEEE products.

As is known to us, the traditional closed-loop supply chain includes one supplier, one manufacturer, one retailer and one recycler. Moreover, the manufacturer is the channel leader in the closed-loop supply chain, while, the retailer is the channel follower. With the development of economy and society, more and more firms participate in the closed-loop supply chain. Thus, the closed-loop supply chain becomes the closed-loop supply chain network which includes multi-supplier, multi-manufacturer and multi-retailer. There is competition relationship between different firms; we think that this competition relationship has advantage to promoting the reuse and recovery of the WEEE products. Based on above theory and practice of the closed-loop supply chain network, with the reward and penalty mechanism and the asymmetric information, we consider five game models which are different in terms of the coordination structure of the reverse supply chain network and power structure of the reward-penalty mechanism: (1) the reverse supply chain network centralized decision-making model; (2) the reverse supply chain network centralized decision-making model with carbon emission constraint; (3) the retailers' competition reverse supply chain network decentralized decision-making model; (4) the retailers' competition reverse supply chain network decentralized decision-making model with carbon emission constraint; and (5) the retailers' competition reverse supply chain network decentralized decision-making model with carbon emission constraint and the government's reward-penalty mechanism. We discuss the following questions:

(1) With the asymmetric information, how does the participate-incentive contract design between the manufacturer and the retailer so that it promotes the retailer to recycling more WEEE products in different decision-making model?

(2) How does the reward and penalty mechanism of the government impact on the recycling rate, the sell price and the profit in different decision-making model?

This paper is different from the other closed-loop supply chain researches in the following aspects. Firstly, we expands the research object which includes one-manufacturer and two-retailer, moreover, there is the competition relationship between two retailers. Assuming that the recycling scale coefficient, the recycling rate and the fixed cost are the retailer one's private information, the manufacturer should identify this private information by designing the participate-incentive contract.

Secondly, the recycling and reuse of the WEEE products can produce much more economy value; the government should participate in this process. We assume that the government as an individual should take part in the process of the recycling of the WEEE products by designing the reward and penalty mechanism. With the asymmetric information, we combine the reward and penalty mechanism with the closed-loop supply chain network and research how the reward and penalty mechanism of the government impact on the sell price and the profit in different decision-making model in the closed-loop supply chain network.

The rest of this paper is organized as follows. In the second section, we present the relevant literature. The theory model and the relevant variables are described in section 3. Five game models which are different in terms of the coordination structure of the reverse supply chain network and power structure of the reward-penalty mechanism are considered in section 4. The result in different decision-making game models can be compared in section 5. Numerical analysis is given in section 6. We can get some conclusions in section 7.

## 2. Literature Review

This paper is closely related to two streams of literature: the coordination strategy with symmetric information when manufacturer or recycler participates in the reverse supply chain; the coordination strategy with symmetric information when the government participates in the reverse supply chain. We briefly review the two streams of related literature.

### 2.1 Coordination strategy with symmetric information when manufacturer or recycler participates in the reverse supply chain

Many scholars discuss the reverse logistics and the closed-loop supply chain. At the beginning of 1990s, some literatures discussed the supply chain coordination problem. For example, Monczka et al.(1998)[1] defined the supply chain coordination problem that it should deal with the relationship between different firms. Mentzer et al.(2001)[2] thought that the supply chain coordination was a union in business operation level. Corbett and Savaskan (2003)[3] pointed out that each member in the reverse supply chain could achieve the coordination of overall interests by making the contract. Wong et al.(2009)[4] researched that it could achieve the coordination by making the contract and improved the efficiency between the JCPenny and Wal-Mart. It could achieve the coordination by many kinds of measures, such as improving the quality of products (Debo et al.,2005[5]), competitive strategy (Majumder and Groenevelt, 2001[6]), technological innovation (Geyer et al. ,2007[7]). Ferrer and Swaminathan (2006, 2010)[8-9] proposed that it should make the different price between the new electronic product and the recycling electronic product.

Ferguson and Toktay (2006)[10] set the nonlinear recovery function and manufacturer could improve their profits by using compulsory acquisition strategy. Assuming that competition between manufacturers, based on the product life cycle theory, Mitra and Webster (2008)[11] analysed the price and coordination strategy in the reverse supply chain. From the aspect of the social responsibility, Yang Yu Xiang et al.(2011)[12] researched the coordination problem in the reverse supply chain. Assuming that manufacturers sell product by using the network market, Zhang Gui Tao et al.(2013)[13] analysed the coordination problem under the case of damaged products. Luo Ding Ti et al.(2010)[14] built the decentralized supply chain consisting of one supplier and two retailers, assuming that they only trade one kinds of product, it could achieve the coordination in supply chain by using the payment incentive measure. Based on the stochastic demand function, Zhang Cui Hua et

al.(2004)[15] built the incentive function under the rewards and punishment mechanism, under the asymmetric information, they built the coordination mechanism between one manufacturer and one supplier. Heese, Cattani, Ferrer et al.(2005)[16] investigated the impact of the recycling waste electronic product on the advantage of oligarch competition in the double oligarch market. Savaskan et al. (2004)[17] analysed efficiency problem in two different types of reverse supply chain. Guide and Teunter (2003)[18] discussed how to set up the demand equilibrium value of waste electronic products, and it could achieve the benefit maximization in the remanufacturing supply chain system. Mukhopadhyay and Setoputro (2005)[19] researched that there was the positive relationship between the performance of electronic products and the best recovery strategy.

**2.2 Coordination strategy with symmetric information when the government participating in reverse supply chain.**

Above literatures mainly discussed how to realize the coordination strategy with manufacturer or recycler participating in the reverse supply chain. The government as a separate entity, and should be involved in the reverse supply chain activities. It was very important that the government promoted the recovery of the WEEE products by making regulations. For example, the government should take carbon tax to the manufacturer for reducing the carbon emission, while the government could make reward and penalty mechanism so that it could promote the manufacturer to improving the recovery and reuse of the WEEE products. Palmer and Walls (1999)[20] pointed out that the government should provide financial and policy to supporting for manufacturer, and encouraging manufacturer to participate in the waste product recycling and remanufacturing activity. Zhang Bao Yin et al. (2006)[21] pointed out that the government should implement certain reward and punishment measures so that the manufacturer could recycle the old product. Dai Yi Sheng et al. (2008)[24] proposed a new fiscal subsidy mechanism, and required the minimum purchase price in the model. Da Qing Li (2008)[23] pointed out that it worthy to be discussed that the government's policy could impact on recycling the waste product.

Wang Wen Bin et al.(2008)[24] discussed the impact of the reward and penalty mechanism on the manufacturer's recycling activitiy. Mitra et al.(2008)[25] proposed a two-stage game model, and manufacturer or recycler existed competition relationship, analysed the government subsidy impact on the reverse supply chain. Based on above analysis, Wang Wen Bin et al.(2009)[26] compared and analysed the difference that the government provided reward and penalty mechanism to manufacturer and collector. Chen et al. (2009)[27] built an environmental regulation pricing strategy, and they concluded that the government should improve the price standard and it could effectively guide manufacturer to improve the recovery rate of waste product. Based on the government's reward and penalty mechanism, Yan Ming et al.(2015)[28] built the WEEE products remanufacturing closed-loop supply chain model, and analysed the impact of the consumer's pay willingness (WTP) for remanufacturing product and government's reward and penalty mechanism on the pricing, recovery rate and the profit. Ne Jia Jia (2015)[29] analysed the variation of the retailer's recovery price, the profit and the total carbon emission in the case of having carbon emission constraints and having not carbon emission constraints. With symmetric information, Cheng Yong Wei et al.(2016)[30] established cooperation decision model for a mixed carbon policy of carbon trading-carbon tax (environmental tax) in a two-stage $S$-$M$ supply chain. With symmetric information, Wang Wen Bin et al. (2016)[31] discussed the impact of the government's reward and penalty mechanism on recycling remanufacturing decision.

The above literatures generally considered the government's reward and penalty mechanism and carbon emission as "instrumental variables" in the reverse supply chain. Moreover, the research object in above literatures included one manufacturer and one retailer. The above literatures considered the impact of the government's reward and penalty mechanism on the manufacturer or the retailer under the symmetric information conditions. On the other hand, this paper analysed the impact of the government's reward and penalty mechanism on the reverse supply chain network under the asymmetric information conditions.

## 3. Theory Model and Variable descriptions

### 3.1. Theory Model

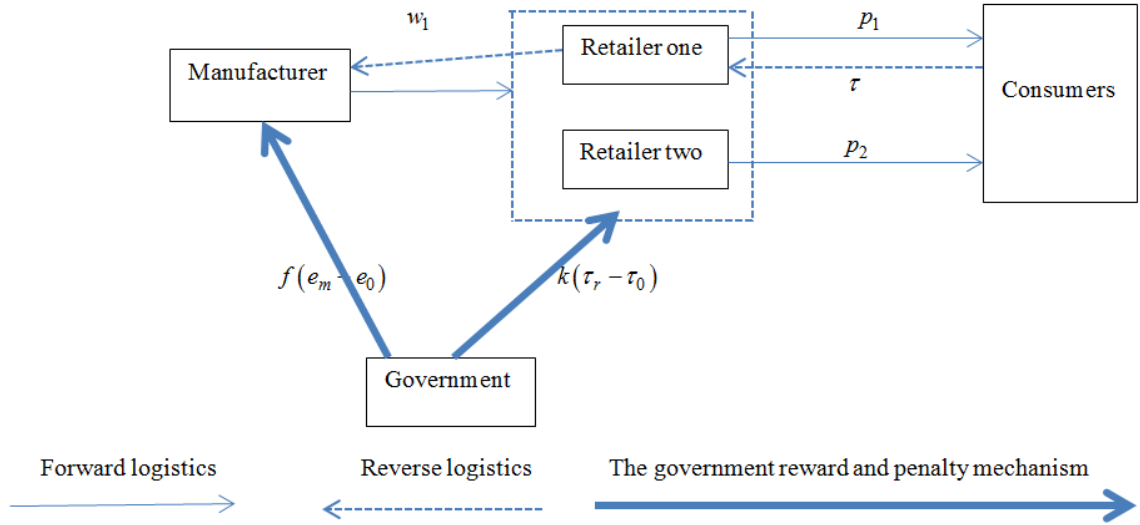

Fig.1. The structure of the reward-penalty mechanism for reverse supply chain network with asymmetric Information and retail-level competition in the presence of carbon emission constraints

In the Fig.1., we consider a single manufacturer (such as computers), two retailers (such as supermarkets), consumers, and government department supply chain network system components. The manufacturer sells product to two retailers. The retailer one recycles the WEEE products, the recycling rate expresses $\tau$, and retailer two doesn't recycle the WEEE products. The manufacturer buy-backs the WEEE products from the retailer one, the buy-back price expresses $w_1$. The manufacturer produces new products by using the old material. The consumer purchases new products, the price can expresses $p_1$ and $p_2$.

The retailer one's fixed cost expresses $I = \beta\tau^2$, $I$ can express the retailer one's fixed cost, it includes recovery network construction funds, advertisement fees and so on, the recycling difficulty coefficient can express $\beta$. The recycling rate can increase with the increasing of the fixed cost. Assuming that the recycling difficulty coefficient, the recycling rate and the fixed cost are the retailer one's private information, it also reflects the asymmetric information between one manufacturer and two retailers. Thus, the recycling difficulty coefficient, the recycling rate and the fixed cost of the retailer one's can express ($H-type$) including $\beta_H, \tau_H, I_H$, and ($L-type$) including $\beta_L, \tau_L, I_L$. Manufacturer doesn't know retailer one's real information about recycling rate, the fixed cost and the recycling difficulty coefficient, however, the manufacturer knows their distribution probability, $P(I_L) = \mu, P(I_H) = 1-\mu$. In order to get more profit, the retailer one may provide false information. Under the asymmetric information conditions, the manufacturer should design the identify contract so that the retailer one can make rational choice.

Assuming that the manufacturer is the channel leader of the reverse supply chain network, the retailer is the channel follower. The government should design the reward and penalty mechanism that the government should set the goal recycling rate and the carbon emission superior limitation. On the one hand, when the recycling rate is lower than the goal recovery rate, retailer one should be punished by the government, on the contrary, the retailer one can get reward. On the other hand, when the carbon emission is larger than the superior limitation, retailer one should be punished by the government, on the contrary, the retailer one can get reward. In this paper, under the asymmetric information conditions, the government provides the reward and penalty mechanism for the carbon emission and recovery rate, when two retailers participate in competition activity, discussing how the government should provide reward and penalty mechanism to promoting two retailers participating in the recycling activity.

### 3.2. Variable Descriptions

In this section, we introduce the meaning of some variables. The retailer one's recycling difficulty coefficient can express $\beta$, $\beta_H$ and $\beta_L$ are the high and low coefficient. The retailer one's recycling rate expresses $\tau$, $\tau_H$ and $\tau_L$ are the high and low recycling rate, which is the retailer one's decision variable. The retailer one's fixed cost can express $I$, $I_H$ and $I_L$ are the high and low fixed cost. The unit testing and sorting cost of manufacturer can express $c_d$. The unit remanufacturing cost of manufacturer can express $c_r$. The manufacturer's unit cost of new product express $c_m$. The retailer one's unit price of new product express $p_1$. The retailer two's unit price of new product is $p_2$. The unit buy-back price of retailer one's can express $w_1$. The market demand of retailer one's can express $q_1$, $q_1 = a - p_1 + \varepsilon p_2$, $a$ is the potential market demand, $\varepsilon$ is the substitute coefficient, $0 < \varepsilon < 1$. The market demand of retailer two's can express $q_2$, $q_2 = a - p_2 + \varepsilon p_1$. The goal recycling rate of the government can express $\tau_0$, when $\tau > \tau_0$, the retailer one can get some reward, and on the contrary, the retailer one can get some penalty. The manufacturer's unit carbon emission can express $e_m$. The goal carbon emission superior limitation can express $e_0$, when the manufacturer's carbon emission is larger than superior limitation, the manufacturer can get penalty, and on the contrary, the manufacturer can get reward. The government unit degree of reward and penalty mechanism for the retailer one can express $f$. The total quantity of the market demand can express $Q$, $Q = q_1 + q_2$. The government total degree of reward and penalty strength for the manufacturer can express $M$, $M = -f[Qe_m - e_0]$. The government unit degree of reward and penalty strength for the manufacturer can express $k$.

## 4. The Decision-Making Model of the Reverse Supply Chain Network Based on the Reward and Penalty Mechanism of the Government under Asymmetric Information

Assuming that the recycling difficulty coefficient, the recycling rate and the fixed cost are the retailer one's private information, it also reflects the asymmetric information between manufacturer and two retailers. The manufacturer should design the identify contract so that the retailer one can make rational choice. The manufacturer design two kinds of identify contracts, includes $\{(w_H, \tau_H), (w_L, \tau_L)\}$. Based on the Principal-agent theory, the manufacturer should set the participation constraint and incentive constraint so that the retailer one can participate in the contract. Five stackelberg game models are built which are different in terms of the coordination structure of the reverse supply chain network and power structure of the reward-penalty mechanism: (1) the reverse supply chain network centralized decision-making model; (2) the reverse supply chain network centralized decision-making model with carbon emission constraint; (3) the retailers' competition reverse supply chain network decentralized decision-making model; (4) the retailers' competition reverse supply chain network decentralized decision-making model with carbon emission constraint; and (5) the retailers' competition reverse supply chain network decentralized decision-making model with carbon emission constraint and the government's reward-penalty mechanism. We discuss the government reward and penalty mechanism design problems in the reverse supply chain network with asymmetric information and retail-level Competition in the presence of carbon emission constraints.

### 4.1 Model Ⅰ: The reverse supply chain network centralized decision-making model

Firstly, we research the reverse supply chain network centralized decision-making model. In this case, in order to identify the true cost of the retailer one, we must use the above contract. In order to comparing with the decentralized decision-making model, assuming that there are one manufacturer and one retailer in the model Ⅰ. The profit of the manufacturer can be expressed as follows:

$$\pi_M^{(1)} = \mu\left\{(a-p_1)\left(p_m - \tau_H^{(1)}\left(w_H^{(1)} + c_d + c_r\right) - \left(1-\tau_H^{(1)}\right)c_m\right)\right\} + (1-\mu)\left\{(a-p_1)\left(p_m - \tau_L^{(1)}\left(w_L^{(1)} + c_d + c_r\right) - \left(1-\tau_L^{(1)}\right)c_m\right)\right\} \quad (1)$$

The profit of the retailer one can be expressed as follows:

$$\pi_{R_1}^{(1)} = \mu\left\{(a-p_1)\tau_H^{(1)}\left(w_H^{(1)} - c\right) + (a-p_1)(p_1 - p_m) - \beta_H\left(\tau_H^{(1)}\right)^2\right\} + (1-\mu)\left\{(a-p_1)\tau_L^{(1)}\left(w_L^{(1)} - c\right) + (a-p_1)(p_1 - p_m) - \beta_L\left(\tau_L^{(1)}\right)^2\right\} \quad (2)$$

The total profit of the reverse supply chain network is that:

$$\pi_T^{(1)} = \mu\left\{(a-p_1)\left(p_1 - \tau_H^{(1)}\left(c_d + c_r - c_m + c\right) - c_m\right) - \beta_H\left(\tau_H^{(1)}\right)^2\right\} + (1-\mu)\left\{(a-p_1)\left(p_1 - \tau_L^{(1)}\left(c_d + c_r - c_m + c\right) - c_m\right) - \beta_L\left(\tau_L^{(1)}\right)^2\right\} \quad (3)$$

The participation constraint and incentive constraint of the retailer one is that:

$$s.t \begin{cases} \tau_L^{(1)}(a-p_1)\left(w_L^{(1)} - c\right) + (a-p_1)(p_1 - p_m) - \beta_L\left(\tau_L^{(1)}\right)^2 \geq \pi_{R_0} \\ \tau_H^{(1)}(a-p_1)\left(w_H^{(1)} - c\right) + (a-p_1)(p_1 - p_m) - \beta_H\left(\tau_H^{(1)}\right)^2 \geq \pi_{R_0} \\ \tau_H^{(1)}(a-p_1)\left(w_H^{(1)} - c\right) + (a-p_1)(p_1 - p_m) - \beta_H\left(\tau_H^{(1)}\right)^2 \geq \tau_L^{(1)}(a-p_1)\left(w_L^{(1)} - c\right) + (a-p_1)(p_1 - p_m) - \beta_H\left(\tau_L^{(1)}\right)^2 \\ \tau_L^{(1)}(a-p_1)\left(w_L^{(1)} - c\right) + (a-p_1)(p_1 - p_m) - \beta_L\left(\tau_L^{(1)}\right)^2 \geq \tau_H^{(1)}(a-p_1)\left(w_H^{(1)} - c\right) + (a-p_1)(p_1 - p_m) - \beta_L\left(\tau_H^{(1)}\right)^2 \end{cases} \quad (4)$$

In equation (4), $\pi_{R_0}$ is the lowest profit of the retailer one, the first two inequalities are the participation constraint, the after two inequalities are the incentive constraint.

The model Ⅰ can be solved by using the Lagrange multiplier method. The optimal solution can express as follows:

$$\tau_H^{(1)} = \frac{\mu(a-p_1)(p_m + c_m - c_r - c - c_d)}{4\beta_H} \quad (5)$$

$$\tau_L^{(1)} = \frac{(1-\mu)(a-p_1)(p_m + c_m - c_r - c - c_d)}{4(\beta_L - \mu\beta_H)} \quad (6)$$

$$w_H^{(1)} = \frac{\mu(p_m + c_m - c_r - c - c_d)}{4} + \frac{(1-\mu)^2(p_m + c_m - c_r - c - c_d)}{4(\beta_L - \mu\beta_H)} + \frac{4\pi_0\beta_H}{(a-p_1)^2(p_m + c_m - c_r - c - c_d)} + c + c_d \quad (7)$$

$$w_L^{(1)} = \frac{(1-\mu)(p_m + c_m - c_r - c - c_d)}{4(\beta_L - \mu\beta_H)} + \frac{(1-\mu)^2(p_m + c_m - c_r - c - c_d)}{4} + \frac{4\pi_0\beta_H(\beta_L - \mu\beta_H)}{(1-\mu)(a-p_1)^2(p_m + c_m - c_r - c - c_d)} + c + c_d \quad (8)$$

### 4.2 Model Ⅱ: The reverse supply chain network centralized decision-making model with carbon emission constraint

In this case, the government total degree of reward and penalty mechanism for the manufacturer is $M = -f[q_1 e_m - e_0]$, in this model, the decision-maker pursue the overall profit maximization. Assuming that there are one manufacturer and one retailer in the model II, the profit of the manufacturer can be expressed as follows:

$$\pi_M^{(2)} = \mu\left\{(a-p_1)\left(p_m - \tau_H^{(2)}\left(w_H^{(1)} + c_d + c_r\right) - \left(1 - \tau_H^{(2)}\right)c_m - fe_m\right) + fe_0\right\} \\ + (1-\mu)\left\{(a-p_1)\left(p_m - \tau_L^{(2)}\left(w_L^{(1)} + c_d + c_r\right) - \left(1 - \tau_L^{(2)}\right)c_m - fe_m\right) + fe_0\right\} \quad (9)$$

The profit of the retailer one can be expressed as follows:

$$\pi_{R_1}^{(2)} = \mu\left\{(a-p_1)\tau_H^{(2)}\left(w_H^{(2)} - c\right) + (a-p_1)(p_1 - p_m) - \beta_H\left(\tau_H^{(2)}\right)^2\right\} + (1-\mu)\left\{(a-p_1)\tau_L^{(2)}\left(w_L^{(2)} - c\right) + (a-p_1)(p_1 - p_m) - \beta_L\left(\tau_L^{(2)}\right)^2\right\} \quad (10)$$

The total profit of the reverse supply chain is that:

$$\pi_T^{(2)} = \mu\left\{(a-p_1)\left(p_1 - \tau_H^{(2)}(c_d + c_r - c_m + c) - c_m - fe_m\right) + fe_0 - \beta_H\left(\tau_H^{(2)}\right)^2\right\} \\ + (1-\mu)\left\{(a-p_1)\left(p_1 - \tau_L^{(2)}(c_d + c_r - c_m + c) - c_m - fe_m\right) + fe_0 - \beta_L\left(\tau_L^{(2)}\right)^2\right\} \quad (11)$$

The participation constraint and incentive constraint of the retailer one is that:

$$s.t \begin{cases} \tau_L^{(2)}(a-p_1)\left(w_L^{(2)} - c\right) + (a-p_1)(p_1 - p_m) - \beta_L\left(\tau_L^{(2)}\right)^2 \geq \pi_{R_0} \\ \tau_H^{(2)}(a-p_1)\left(w_H^{(2)} - c\right) + (a-p_1)(p_1 - p_m) - \beta_H\left(\tau_H^{(2)}\right)^2 \geq \pi_{R_0} \\ \tau_H^{(2)}(a-p_1)\left(w_H^{(2)} - c\right) + (a-p_1)(p_1 - p_m) - \beta_H\left(\tau_H^{(2)}\right)^2 \geq \tau_L^{(2)}(a-p_1)\left(w_L^{(2)} - c\right) + (a-p_1)(p_1 - p_m) - \beta_H\left(\tau_L^{(2)}\right)^2 \\ \tau_L^{(2)}(a-p_1)\left(w_L^{(2)} - c\right) + (a-p_1)(p_1 - p_m) - \beta_L\left(\tau_L^{(2)}\right)^2 \geq \tau_H^{(2)}(a-p_1)\left(w_H^{(2)} - c\right) + (a-p_1)(p_1 - p_m) - \beta_L\left(\tau_H^{(2)}\right)^2 \end{cases} \quad (12)$$

In equation (12), $\pi_{R_0}$ is the lowest profit of the retailer one, the first two inequalities is the participation constraint, the after two inequalities is the incentive constraint.

The model II can be solved by using the Lagrange multiplier method. The optimal solution can express as follows:

$$\tau_H^{(2)} = \frac{\mu(a-p_1)(p_m + c_m - c_r - c - c_d) - fe_m}{4\beta_H} \quad (13)$$

$$\tau_L^{(2)} = \frac{(1-\mu)(a-p_1)(p_m + c_m - c_r - c - c_d) - fe_m}{4(\beta_L - \mu\beta_H)} \quad (14)$$

$$w_H^{(2)} = \frac{\mu(p_m + c_m - c_r - c - c_d) - fe_m}{4} + \frac{(1-\mu)^2(p_m + c_m - c_r - c - c_d) + (\beta_L - \mu\beta_H - fe_m)}{4(\beta_L - \mu\beta_H)} + \frac{4\pi_0\beta_H}{(a-p_1)^2(p_m + c_m - c_r - c - c_d)} + c + c_d \quad (15)$$

$$w_L^{(2)} = \frac{(1-\mu)(p_m + c_m - c_r - c - c_d) - fe_m}{4(\beta_L - \mu\beta_H)} + \frac{(1-\mu)^2(p_m + c_m - c_r - c - c_d)}{4} + \frac{4\pi_0\beta_H(\beta_L - \mu\beta_H)}{(1-\mu)(a-p_1)^2(p_m + c_m - c_r - c - c_d)} + c + c_d \quad (16)$$

### 4.3 Model III: the retailers' competition reverse supply chain network decentralized decision-making model

In this case, the manufacturer is the channel leader in the Stackelberg game model, while, two retailers are the follower in the Stackelberg game model. The manufacturer commissioned the retailer one recycling the WEEE products. At the same time, there is competition relationship between two retailers. In model III, the manufacturer and two retailers should consider their own maximize profit. The profit of the manufacturer can be expressed as follows:

$$\pi_M^{(3)} = \mu\left\{\left(2a-(p_1+p_2)+\varepsilon(p_1+p_2)\right)\left(p_m - \tau_H^{(3)}\left(w_H^{(3)}+c_d+c_r\right)-\left(1-\tau_H^{(3)}\right)c_m\right)\right\} \quad (17)$$
$$+(1-\mu)\left\{\left(2a-(p_1+p_2)+\varepsilon(p_1+p_2)\right)\left(p_m - \tau_L^{(3)}\left(w_L^{(3)}+c_d+c_r\right)-\left(1-\tau_L^{(3)}\right)c_m\right)\right\}$$

The participation and incentive constraint of the retailer one is that:

$$s.t \begin{cases} \tau_L^{(3)}(a-p_1+\varepsilon p_2)\left(w_L^{(3)}-c\right)+(a-p_1+\varepsilon p_2)(p_1-p_m)-\beta_L\left(\tau_L^{(3)}\right)^2 \geq \pi_{R_0} \\ \tau_H^{(3)}(a-p_1+\varepsilon p_2)\left(w_H^{(3)}-c\right)+(a-p_1+\varepsilon p_2)(p_1-p_m)-\beta_H\left(\tau_H^{(3)}\right)^2 \geq \pi_{R_0} \\ \tau_H^{(3)}(a-p_1+\varepsilon p_2)\left(w_H^{(3)}-c\right)+(a-p_1+\varepsilon p_2)(p_1-p_m)-\beta_H\left(\tau_H^{(3)}\right)^2 \geq \tau_L^{(3)}(a-p_1+\varepsilon p_2)\left(w_L^{(3)}-c\right)+(a-p_1+\varepsilon p_2)(p_1-p_m)-\beta_H\left(\tau_L^{(3)}\right)^2 \\ \tau_L^{(3)}(a-p+\varepsilon p_2)\left(w_L^{(3)}-c\right)+(a-p_1+\varepsilon p_2)(p_1-p_m)-\beta_L\left(\tau_L^{(3)}\right)^2 \geq \tau_H^{(3)}(a-p_1+\varepsilon p_2)\left(w_H^{(3)}-c\right)+(a-p_1+\varepsilon p_2)(p_1-p_m)-\beta_L\left(\tau_H^{(3)}\right)^2 \end{cases} \quad (18)$$

The retailer one is responsible for recycling the WEEE product, the decision variables are $\tau$ and $p_1$, the retailer one's profit can be expressed as follows:

$$\pi_{r_1}^{(3)} = \mu\left\{(a-p_1+\varepsilon p_2)\tau_H^{(3)}\left(w_H^{(3)}-c\right)+(a-p_1+\varepsilon p_2)(p_1-p_m)-\beta_H\left(\tau_H^{(3)}\right)^2\right\} \quad (19)$$
$$+(1-\mu)\left\{(a-p_1+\varepsilon p_2)\tau_L^{(3)}\left(w_L^{(3)}-c\right)+(a-p_1+\varepsilon p_2)(p_1-p_m)-\beta_L\left(\tau_L^{(3)}\right)^2\right\}$$

The retailer two is not responsible for recycling the WEEE product, the decision variable is $p_2$, the retailer two's profit can be expressed as follows:

$$\pi_{r_2}^{(3)} = (p_2-p_m)(a-p_2+\varepsilon p_1) \quad (20)$$

The decision order for the game model III: Firstly, the manufacturer should decide the buy-back price $w_H^{(3)}$ and $w_L^{(3)}$, then, the retailer one should decide the recovery rate $\tau_H^{(3)}$ and $\tau_L^{(3)}$, lastly, two retailers decide their prices $p_{R_1}^{(3)}$ and $p_{R_2}^{(3)}$.

The model III can be solved by using the Lagrange multiplier method. The optimal solution can express as follows:

$$w_H^{(3)} = \frac{\beta_L\left(2a+2c-2p_m+a\varepsilon-c\varepsilon p_m-c\varepsilon^2 p_m+\varepsilon^2 p_m-ac\varepsilon-2ac\right)}{2\beta_H\left(4\beta_L-c^2-\beta_L^2+\beta_L\mu\right)-2\beta_L\left(c^2\mu-\mu\beta_H^2\right)} \quad (21)$$

$$w_L^{(3)} = \frac{\beta_H\left(2a\beta_L-c\varepsilon p_m+\varepsilon\beta_L p_m-c\varepsilon^2 p_m+\varepsilon^2\beta_L p_m-ac\varepsilon+a\varepsilon\beta_L-2ac\right)+2\beta_L(cp_m-\beta_L p_m)}{2\beta_H\left(\beta_L-c^2-\beta_L^2+\beta_L\mu\right)-2\beta_L\left(c^2\mu-\mu\beta_H^2\right)} \quad (22)$$

$$\tau_H^{(3)} = \frac{\pi_{R_0}\mu(a-\varepsilon)\left(2a-2c+a\varepsilon+c\varepsilon^2\right)}{\mu\beta_H(\varepsilon-2\mu)^2} \quad (23)$$

$$\tau_L^{(3)} = \frac{\pi_{R_0}(1-\mu)\left(2a-2c+a\varepsilon+c\varepsilon^2\right)}{\mu(\beta_H-\beta_L)(\varepsilon^2-4\mu)} \quad (24)$$

$$p_{R_1}^{(3)} = \frac{\mu(2+\varepsilon)\left[2\beta_H-(p_m+c_m+c_r+c_d-c)^2\right]+(1-\mu)(4+\varepsilon)\left[2\beta_L-(p_m+c_m+c_r+c_d-c)^2\right]}{4\beta_H-(p_m+c_m+c_r+c_d-c)^2+4\beta_L-\varepsilon^2} \quad (25)$$

$$p_{R_2}^{(3)} = \frac{\mu(6+\varepsilon)\left[4(\beta_H)^2-2(p_m+c_m+c_r-c)^2\right]+(1-\mu)(8+2\varepsilon)\left[4(\beta_L)^2-2(p_m+c_m+c_r-c)^2\right]}{4\beta_H-(p_m+c_m+c_r+c_d-c)^2+4\beta_L-\varepsilon^2} \quad (26)$$

### 4.4 Model IV: the retailers' competition reverse supply chain network decentralized decision-making model with carbon emission constraint

In this case, the manufacturer is the channel leader in the Stackelberg game model, while, two retailers are the channel followers. The manufacturer commissioned the retailer one recycling the WEEE products. At the same time, there is competition relationship between two retailers'. The government total degree of reward and penalty strength to the manufacturer is $M = -f[q_1 e_m - e_0]$, the profit of the manufacturer can be expressed as follows

$$\pi_M^{(4)} = \mu\left\{(2a-(p_1+p_2)+\varepsilon(p_1+p_2))\left(p_m - \tau_H^{(4)}\left(w_H^{(4)} + c_d + c_r\right) - \left(1-\tau_H^{(4)}\right)c_m - fe_m\right) + fe_0\right\} \\ + (1-\mu)\left\{(2a-(p_1+p_2)+\varepsilon(p_1+p_2))\left(p_m - \tau_L^{(4)}\left(w_L^{(4)} + c_d + c_r\right) - \left(1-\tau_L^{(4)}\right)c_m - fe_m\right) + fe_0\right\} \tag{27}$$

The participation and incentive constraint of the retailer one is that:

$$s.t \begin{cases} \tau_L^{(4)}(a-p_1+\varepsilon p_2)\left(w_L^{(4)}-c\right)+(a-p_1+\varepsilon p_2)(p_1-p_m)-\beta_L\left(\tau_L^{(4)}\right)^2 \geq \pi_{R_0} \\ \tau_H^{(4)}(a-p_1+\varepsilon p_2)\left(w_H^{(4)}-c\right)+(a-p_1+\varepsilon p_2)(p_1-p_m)-\beta_H\left(\tau_H^{(4)}\right)^2 \geq \pi_{R_0} \\ \tau_H^{(4)}(a-p_1+\varepsilon p_2)\left(w_H^{(4)}-c\right)+(a-p_1+\varepsilon p_2)(p_1-p_m)-\beta_H\left(\tau_H^{(4)}\right)^2 \geq \tau_L^{(4)}(a-p_1+\varepsilon p_2)\left(w_L^{(4)}-c\right)+(a-p_1+\varepsilon p_2)(p_1-p_m)-\beta_H\left(\tau_L^{(4)}\right)^2 \\ \tau_L^{(4)}(a-p_1+\varepsilon p_2)\left(w_L^{(4)}-c\right)+(a-p_1+\varepsilon p_2)(p_1-p_m)-\beta_L\left(\tau_L^{(4)}\right)^2 \geq \tau_H^{(4)}(a-p_1+\varepsilon p_2)\left(w_H^{(4)}-c\right)+(a-p_1+\varepsilon p_2)(p_1-p_m)-\beta_L\left(\tau_H^{(4)}\right)^2 \end{cases} \tag{28}$$

The retailer one's profit can be expressed as follows:

$$\pi_{R_1}^{(4)} = \mu\left\{(a-p_1+\varepsilon p_2)\tau_{R_1^H}^{(4)}\left(w_{M^H}^{(4)}-c\right)+(a-p_1+\varepsilon p_2)(p_1-p_m)-\beta_H\left(\tau_{R_1^H}^{(4)}\right)^2\right\} \\ +(1-\mu)\left\{(a-p_1+\varepsilon p_2)\tau_{R_1^L}^{(4)}\left(w_{M^L}^{(4)}-c\right)+(a-p_1+\varepsilon p_2)(p_1-p_m)-\beta_L\left(\tau_{R_1^L}^{(4)}\right)^2\right\} \tag{29}$$

The retailer two's profit can be expressed as follows:

$$\pi_{R_2}^{(4)} = (p_2 - p_m)(a - p_2 + \varepsilon p_1) \tag{30}$$

The model IV can be solved by using the Lagrange multiplier method. The optimal solution can express as follows:

$$w_H^{(4)} = \frac{\beta_L\left(2a+2c-2p_m+a\varepsilon-c\varepsilon p_m-c\varepsilon^2 p_m+\varepsilon^2 p_m-ac\varepsilon-2ac+fe_m\right)-fe_0}{2\beta_H\left(4\beta_L-c^2-\beta_L^2+\beta_L\mu\right)-2\beta_L\mu\left(c^2-\beta_H^2\right)} \tag{31}$$

$$w_L^{(4)} = \frac{\beta_H\left(2a\beta_L-c\varepsilon p_m+\varepsilon\beta_L p_m-c\varepsilon^2 p_m+\varepsilon^2\beta_L p_m-ac\varepsilon+a\varepsilon\beta_L-2ac\right)+\beta_L\left(2cp_m-2\beta_L p_m+fe_m\right)-fe_0}{2\beta_H\left(4\beta_L-c^2-\beta_L^2+\beta_L\mu\right)-2\beta_L\mu\left(c^2-\beta_H^2\right)} \tag{32}$$

$$\tau_{R_1^H}^{(4)} = \frac{\pi_{R_0}\mu(a-\varepsilon)\left(2a-2c+a\varepsilon+c\varepsilon^2\right)+\beta_H fe_m}{\mu\beta_H(\varepsilon-2\mu)^2} \tag{33}$$

$$\tau_{R_1^L}^{(4)} = \frac{\pi_{R_0}(1-\mu)\left(2a-2c+a\varepsilon+c\varepsilon^2\right)+(\beta_H-\beta_L)fe_m}{\mu(\beta_H-\beta_L)(\varepsilon^2-4\mu)} \tag{34}$$

$$p_{R_1}^{(4)} = \frac{\mu(2+\varepsilon)(a-\varepsilon)\left[2\beta_H+(p_m+c_m+c_r+c_d-c)^2\right]+(a+\varepsilon)(1-\mu)(4+\varepsilon)\left[2\beta_L-(p_m+c_m+c_r+c_d-c)^2\right]}{4\beta_H-(a-\varepsilon)(p_m+c_m+c_r+c_d-c)^2+4\beta_L-\varepsilon^2} \\ + \frac{fe_m}{4\beta_H+(a-\varepsilon)(p_m+c_m+c_r+c_d-c)^2} + fe_0 \tag{35}$$

$$p_{R_2}^{(4)} = \frac{\mu(a-\varepsilon)(6+\varepsilon)\left[4(\beta_H)^2-2(p_m+c_m+c_r-c)^2\right]+(a+\varepsilon)(1-\mu)(8+2\varepsilon)\left[4(\beta_L)^2-2(p_m+c_m+c_r-c)^2\right]}{4\beta_H-(a-\varepsilon)(p_m+c_m+c_r+c_d-c)^2+4\beta_L-\varepsilon^2} \\ + \frac{fe_m}{(a-\varepsilon)(p_m+c_m+c_r+c_d-c)^2} + \frac{fe_0}{\beta_H(a-\varepsilon)} \tag{36}$$

**4.5 Model V: the retailers' competition reverse supply chain network decentralized decision-making model with carbon emission constraint and the government's reward-penalty mechanism**

In this case, the government should give carbon emission constraints for the manufacturer, moreover, the government should give reward and penalty mechanism for the retailer one's recycling rate. Thus, the retailer one can recycle more the WEEE products. The total degree of the reward and penalty strength is $k(\tau_r - \tau_0)$.

The retailer two does not recycle the WEEE products, the retailer two can get penalty by the government, and the total degree of the penalty strength is $k\tau_0$. The profit of the manufacturer can be expressed as follows:

$$\pi_M^{(5)} = \mu\{(2a-(p_1+p_2)+\varepsilon(p_1+p_2))(p_m - \tau_H(w_H + c_d + c_r) - (1-\tau_H)c_m - fe_m) + fe_0\} \\ + (1-\mu)\{(2a-(p_1+p_2)+\varepsilon(p_1+p_2))(p_m - \tau_L(w_L + c_d + c_r) - (1-\tau_L)c_m - fe_m) + fe_0\} \quad (37)$$

The participation and incentive constraint of the retailer one is that:

$$s.t \begin{cases} \tau_L^{(5)}(a-p_1+\varepsilon p_2)(w_L^{(5)}-c)+(a-p_1+\varepsilon p_2)(p_1-p_m)-\beta_L\left(\tau_L^{(5)}\right)^2+k(\tau_r-\tau_0) \geq \pi_{R_0} \\ \tau_H^{(5)}(a-p_1+\varepsilon p_2)(w_H^{(5)}-c)+(a-p_1+\varepsilon p_2)(p_1-p_m)-\beta_H\left(\tau_H^{(5)}\right)^2+k(\tau_r-\tau_0) \geq \pi_{R_0} \\ \tau_H^{(5)}(a-p_1+\varepsilon p_2)(w_H^{(5)}-c)+(a-p_1+\varepsilon p_2)(p_1-p_m)-\beta_H\left(\tau_H^{(5)}\right)^2+k(\tau_r-\tau_0) \geq \tau_L^{(5)}(a-p_1+\varepsilon p_2)(w_L^{(5)}-c)+(a-p_1+\varepsilon p_2)(p_1-p_m)-\beta_H\left(\tau_L^{(5)}\right)^2+k(\tau_r-\tau_0) \\ \tau_L^{(5)}(a-p_1+\varepsilon p_2)(w_L^{(5)}-c)+(a-p_1+\varepsilon p_2)(p_1-p_m)-\beta_L\left(\tau_L^{(5)}\right)^2+k(\tau_r-\tau_0) \geq \tau_H^{(5)}(a-p_1+\varepsilon p_2)(w_H^{(5)}-c)+(a-p_1+\varepsilon p_2)(p_1-p_m)-\beta_L\left(\tau_H^{(5)}\right)^2+k(\tau_r-\tau_0) \end{cases} \quad (38)$$

The retailer one's profit can be expressed as follows:

$$\pi_{R_1}^{(5)} = \mu\left\{(a-p_1+\varepsilon p_2)\tau_H^{(5)}\left(w_H^{(5)}-c\right)+(a-p_1+\varepsilon p_2)(p_1-p_m)-\beta_H\left(\tau_H^{(5)}\right)^2+k(\tau_r-\tau_0)\right\} \\ +(1-\mu)\left\{(a-p_1+\varepsilon p_2)\tau_L^{(5)}\left(w_L^{(5)}-c\right)+(a-p_1+\varepsilon p_2)(p_1-p_m)-\beta_L\left(\tau_L^{(5)}\right)^2+k(\tau_r-\tau_0)\right\} \quad (39)$$

The retailer two's profit can be expressed as follows:

$$\pi_{R_2}^{(5)} = (p_2 - p_m)(a - p_2 + \varepsilon p_1) - k\tau_0 \quad (40)$$

The model V can be solved by using the Lagrange multiplier method. The optimal solution can express as follows:

$$w_H^{(5)} = \frac{\beta_L(2a+2c-2p_m+a\varepsilon-c\varepsilon p_m-c\varepsilon^2 p_m+\varepsilon^2 p_m-a c\varepsilon-2ac+fe_m)-fe_0}{2\beta_H(4\beta_L-c^2-\beta_L^2+\beta_L\mu)-2\mu\beta_L(c^2-\beta_H^2)} \\ + \frac{k(2+\varepsilon^2-\mu^2)}{2\beta_H(4\beta_L-c^2-w_L^2)-2c^2\beta_L\mu} \quad (41)$$

$$w_L^{(5)} = \frac{\beta_H(2a\beta_L-c\varepsilon p_m+\varepsilon\beta_L p_m-c\varepsilon^2 p_m+\varepsilon^2\beta_L p_m-ac\varepsilon-2ac+a\varepsilon\beta_L)+\beta_L(2cp_m-2\beta_L p_m+fe_m)-fe_0}{2\beta_H(4\beta_L-c^2-\beta_L^2+\beta_L\mu)-2\beta_L\mu(c^2-\beta_H^2)} \\ + \frac{k(2-\varepsilon^2-2\mu^2)}{2\beta_H(4\beta_L-c^2-\beta_L^2+\beta_L\mu)} \quad (42)$$

$$\tau_H^{(5)} = \frac{\pi_{R_0}\mu(a-\varepsilon)(2a-2c+a\varepsilon+c\varepsilon^2)+\beta_H fe_m}{\mu\beta_H(\varepsilon-2\mu)^2} + \frac{(8-\varepsilon^2)k}{(\varepsilon-2\mu)^2} \quad (43)$$

$$\tau_L^{(5)} = \frac{\pi_{R_0}(1-\mu)(2a-2c+a\varepsilon+c\varepsilon^2)+(\beta_H-\beta_L)fe_m}{\mu(\beta_H-\beta_L)(\varepsilon^2-4\mu)} + \frac{(4-\varepsilon^2)k-\varepsilon}{(\varepsilon^2-4\mu)} \quad (44)$$

$$p_{R_1}^{(5)} = \frac{\mu(2+\varepsilon)(a-\varepsilon)\left[2\beta_H+(p_m+c_m+c_r+c_d-c)^2\right]+(a+\varepsilon)(1-\mu)(4+\varepsilon)\left[2\beta_L-(p_m+c_m+c_r+c_d-c)^2\right]}{4\beta_H-(a-\varepsilon)(p_m+c_m+c_r+c_d-c)^2+4\beta_L-\varepsilon^2} \\ + \frac{fe_m}{4\beta_H+(a-\varepsilon)(p_m+c_m+c_r+c_d-c)^2} + \frac{2fe_0-(c_m+c_r+c_d-c)k}{(p_m+c_m+c_r+c_d-c)^2} \quad (45)$$

$$p_{R_2}^{(5)} = \frac{\mu(a-\varepsilon)(6+\varepsilon)\left[4(\beta_H)^2-2(p_m+c_m+c_r-c)^2\right]+(a+\varepsilon)(1-\mu)(8+2\varepsilon)\left[4(\beta_L)^2-2(p_m+c_m+c_r-c)^2\right]}{4\beta_H-(a-\varepsilon)(p_m+c_m+c_r+c_d-c)^2+4\beta_L-\varepsilon^2} \\ + \frac{fe_m}{(a-\varepsilon)(p_m+c_m+c_r+c_d-c)^2} + \frac{3fe_0-(c_m+c_r+c_d-c)(k-\varepsilon)}{(p_m+c_m+c_r+c_d-c)^2} \quad (46)$$

## 5. Comparing the result in different decision-making game models

In this section, we mainly compare the result in different decision-making game models; we can get the following propositions:

**Proposition 1.** Comparing the model I with the model II: If $fe_m < (\beta_L - \mu\beta_H)/(1+\beta_L - \mu\beta_H)$, thus $w_H^{(2)} - w_H^{(1)} > 0$; if $fe_m > (\beta_L - \mu\beta_H)/(1+\beta_L - \mu\beta_H)$, thus $w_H^{(2)} - w_H^{(1)} < 0$; Because $-fe_m < 0$, if $\beta_L - \mu\beta_H > 0$, thus $w_L^{(2)} - w_L^{(1)} < 0$.

Proposition 1 suggests that the retailer one's fixed cost is higher, if $fe_m < (\beta_L - \mu\beta_H)/(1+\beta_L - \mu\beta_H)$, in other words, when the manufacturer's carbon emission is smaller than the total degree of the reward and penalty strength, the buy-back price in the model II is larger than that in the model I. On the contrary, if $fe_m > (\beta_L - \mu\beta_H)/(1+\beta_L - \mu\beta_H)$, in other words, when the manufacturer's carbon emission is larger than the total degree of the reward and penalty strength, the buy-back price in the model II is smaller than that in the model I. The retailer one's fixed cost is lower, the buy-back price in the model II is lower than that in the model I.

**Proposition 2.** Comparing the model I with the model II: If $\beta_H > 0$, $-fe_m < 0$, thus $\tau_H^{(2)} - \tau_H^{(1)} < 0$; if $-fe_m < 0$, $\beta_L - \mu\beta_H > 0$, thus $\tau_L^{(2)} - \tau_L^{(1)} < 0$.

Proposition 2 suggests that no matter how high or low of the fixed cost in the retailer one, the recycling rate in the model II is smaller than that in the model I.

**Proposition 3.** Comparing the model III with the model IV: No matter how high or low of the fixed cost in the retailer one, if it satisfy the following condition ① $e_m < \frac{e_0}{\beta_L}$, ② $\beta_H(4\beta_L - c^2 - \beta_L^2 + \mu\beta_L) < \beta_L(\mu c^2 - \mu\beta_H^2)$, the buy-back price in the model IV is larger than that in the model III. On the contrary, if it satisfy the following condition ① $e_m > \frac{e_0}{\beta_L}$, ② $\beta_H(4\beta_L - c^2 - \beta_L^2 + \mu\beta_L) < \beta_L(\mu c^2 - \mu\beta_H^2)$, the buy-back price in the model IV is smaller than that in the model III. Proposition 3 suggest that the retailer one's buy-back price decide on the manufacturer's carbon emission $e_m$ when the government implements the reward and penalty mechanism for the manufacturer's carbon emission.

**Proposition 4.** Comparing the model IV with the model III: If $\beta_H fe_m > 0, \mu > 0, \beta_H > 0, (\varepsilon - 2\mu)^2 > 0$, thus $\tau_H^{(4)} - \tau_H^{(3)} > 0$; If $\beta_H > \beta_L, \mu > 0$, when $\varepsilon^2 - 4\mu > 0$, thus $\tau_L^{(4)} - \tau_L^{(3)} > 0$, when $\varepsilon^2 - 4\mu < 0$, thus $\tau_L^{(4)} - \tau_L^{(3)} < 0$.

Proposition 4 suggest that the retailer one's fixed cost is high, the recycling rate in the model IV is larger than that in the model III. When the retailer one's fixed cost is low, if $\varepsilon^2 - 4\mu > 0$, the recycling rate in the model IV is larger than that in the model III, if $\varepsilon^2 - 4\mu < 0$, the recycling rate in the model IV is lower than that in the model III.

**Proposition 5.** In the model IV, if $fe_m > 0, fe_0 > 0, a - \varepsilon > 0, (p_m + c_m + c_r + c_d - c)^2 > 0$, thus, $p_{R_1}^{(4)} - p_{R_1}^{(3)} > 0$, $p_{R_2}^{(4)} - p_{R_2}^{(3)} > 0$. Proposition 5 suggest that it can improve two retailer's price that the government implement carbon emissions constraints for the manufacturer.

Proving that $p_{R_1}^{(4)} - p_{R_1}^{(3)} = \dfrac{fe_m}{4\beta_H + (a-\varepsilon)(p_m+c_m+c_r+c_d-c)^2} + fe_0$, because $fe_0 > 0, a-\varepsilon > 0, (p_m+c_m+c_r+c_d-c)^2 > 0$

$p_{R_2}^{(4)} - p_{R_2}^{(3)} = \dfrac{fe_m}{(a-\varepsilon)(p_m+c_m+c_r+c_d-c)^2} + \dfrac{fe_0}{\beta_H(a-\varepsilon)}$, if $fe_m > 0, fe_0 > 0, a-\varepsilon > 0, (p_m+c_m+c_r+c_d-c)^2 > 0$, thus $p_{R_1}^{(4)} - p_{R_1}^{(3)} > 0$, $p_{R_2}^{(4)} - p_{R_2}^{(3)} > 0$.

**Proposition 6.** Comparing the model Ⅴ with the model Ⅳ: If the retailer one's fixed cost is high, when it satisfy the following condition ① $4\beta_L - c^2 - \beta_L^2 > 0$ ② $(4\beta_L - c^2 - \beta_L^2) > c^2\mu$, the buy-back price in the model Ⅴ is larger than that in the model Ⅳ. If the retailer one's fixed cost is low, when if satisfy the following condition ① $2 - \varepsilon^2 - 2\mu^2 > 0$ ② $4\beta_L - c^2 - \beta_L^2 + \beta_L\mu > 0$, the buy-back price in the model Ⅴ is larger than that in the model Ⅳ.

**Proposition 7.** Comparing the model Ⅴ with the model Ⅳ: If $8 - \varepsilon^2 > 0, k > 0, (\varepsilon - 2\mu)^2 > 0$, thus, $\tau_H^{(5)} - \tau_H^{(4)} > 0$. If $\varepsilon^2 - 4\mu > 0, k > \varepsilon^2/4 - \varepsilon^2$, thus $\tau_L^{(5)} - \tau_L^{(4)} > 0$. If $\varepsilon^2 - 4\mu > 0$, $k < \varepsilon^2/4 - \varepsilon^2$, thus $\tau_L^{(5)} - \tau_L^{(4)} < 0$.

Proposition 7 suggest that the retailer one's fixed cost is high, the recycling rate in the model Ⅴ is larger than that in the model Ⅳ. The retailer one's fixed cost is low, if it satisfies the condition $\varepsilon^2 - 4\mu > 0$ and $k > \dfrac{\varepsilon^2}{4-\varepsilon^2}$, the recycling rate in the model Ⅴ is larger than that in the model Ⅳ. On the contrary, if it satisfies the condition $\varepsilon^2 - 4\mu > 0$ and $k < \dfrac{\varepsilon^2}{4-\varepsilon^2}$, the recycling rate in the model Ⅴ is smaller than that in the model Ⅳ.

**Proposition 8.** Comparing the model Ⅴ with the model Ⅳ: The retailer one's price is that $k < \dfrac{[2-(p_m+c_m+c_r+c_d-c)^2]fe_0}{c_m+c_r+c_d-c}$, the retailer two's price is that $k < \dfrac{3fe_0}{c_m+c_r+c_d-c} - \dfrac{fe_0(p_m+c_m+c_r+c_d-c)^2}{\beta_H(a-\varepsilon)(c_m+c_r+c_d-c)} + \varepsilon$. Proposition 8 suggest that the two retailers' price $p_1$ and $p_2$ can be impacted by the $k$. No matter how large or low the fixed cost of the retailer one is, Only if the unit degree of the reward and penalty strength reaches certain level, the retailer's price in the model Ⅴ is lower than that in the model Ⅳ. On the contrary, the retailer's price in the model Ⅴ is larger than that in the model Ⅳ.

## 6. Numerical analysis

In above section, we analysis that with the asymmetric information, how does the participate-incentive contract design between the manufacturer and the retailer so that it can promote the retailer to recycling more WEEE products in different decision-making model? How does the reward and penalty mechanism of the government impact on the recycling rate, the sell price and the profit in different decision-making model? Because the decision-making model is very complex in the model Ⅳ and Ⅴ, in this section, we analysis the impact of the government reward and penalty mechanism on the recycling rate, the buy-back price and the retailer's price. The relevant parameter is set as follows:

$c = 4$, $\beta_H = 0.7$, $\beta_L = 0.5$, $I_H = 40$, $I_L = 30$, $c_d = 3$, $c_r = 2.6$, $c_m = 2$, $p_1 = 1.7$, $p_2 = 1.9$, $a = 3$, $\varepsilon = 0.4$, $\tau_0 = 0.8$, $e_m = 0.9$, $e_0 = 1.3$, $p_m = 1.3$, $f \in [3,10], k \in [2,8]$. The decision result in the model Ⅳ and Ⅴ can be expressed in the table 1.

**Table 1 the decision result in the model Ⅳ and Ⅴ**

| $f$ | $k$ | $w_{H_4}^{M*}$ | $w_{L_4}^{M*}$ | $\tau_{H_4}^{R*}$ | $\tau_{L_4}^{R*}$ | $w_{H_5}^{M*}$ | $w_{L_5}^{M*}$ | $\tau_{H_5}^{R*}$ | $\tau_{L_5}^{R*}$ |
|---|---|---|---|---|---|---|---|---|---|
| 3 | 2 | 3.8 | 2.9 | 0.38 | 0.21 | 4.5 | 3.3 | 0.56 | 0.43 |
| 5 | 4 | 4.0 | 3.2 | 0.42 | 0.24 | 5.3 | 3.7 | 0.63 | 0.47 |
| 7 | 6 | 4.1 | 3.3 | 0.45 | 0.28 | 5.9 | 4.6 | 0.69 | 0.56 |
| 9 | 8 | 4.3 | 3.5 | 0.49 | 0.32 | 6.8 | 5.5 | 0.78 | 0.63 |

We can get the following conclusions:

(1) When the retailer one's fixed cost is high, in the model Ⅴ, the buy-back price of the manufacturer's can increase with the increasing of the degree of the reward and penalty strength, the range between 4 and 8. On the other hand, in the model Ⅳ, the buy-back price of the manufacturer's can increase with the increasing of the degree of the reward and penalty strength, the range between 3 and 5. When the retailer one's fixed cost is low, in the model Ⅴ, the buy-back price of the manufacturer's can increase with the increasing of the degree of the reward and penalty strength, the range between 3 and 6. On the other hand, in the model Ⅳ, the buy-back price of the manufacturer's can increase with the increasing of the degree of the reward and penalty strength, the range between 3 and 4. This suggests that it can promote the increasing of the buy-back price by implementing the reward and penalty mechanism.

(2) When the retailer one's fixed cost is high, in the model Ⅴ, the recycling rate of the retailer one can increase with the increasing of the degree of the reward and penalty strength, the range between 0.5 and 0.8. On the other hand, in the model Ⅳ, the recycling rate of the retailer one can increase with the increasing of the degree of the reward and penalty strength, the range between 0.35 and 0.5. When the retailer one's fixed cost is low, in the model Ⅴ, the recycling rate of the retailer one can increase with the increasing of the degree of the reward and penalty strength, the range between 0.4 and 0.7. On the other hand, in the model Ⅳ, the recycling rate of the retailer one can increase with the increasing of the degree of the reward and penalty strength, the range between 0.2 and 0.35. This suggests that it can promote the increasing of the recycling rate by implementing the reward and penalty mechanism.

## 7. Conclusion

In this paper, we discuss the government's reward and penalty mechanism in the presence of asymmetric information and carbon emission constraint when downstream retailers compete in a reverse supply chain network. Considering five game models which are different in terms of the coordination structure of the reverse supply chain network and power structure of the reward-penalty mechanism: (1) the reverse supply chain network centralized decision-making model; (2) the reverse supply chain network centralized decision-making model with carbon emission constraint; (3) the retailers' competition reverse supply chain network decentralized decision-making model; (4) the retailers' competition reverse supply chain network decentralized decision-making model with carbon emission constraint; (5) the retailers' competition reverse supply chain network decentralized decision-making model with carbon emission constraint and the government's reward-penalty mechanism. Building the participation-incentive contract under each model use the principal-agent theory and solving the model use the Lagrange multiplier method. We can get the following conclusion: (1) when the government implements the reward-penalty mechanism for carbon emission and recycling simultaneously, the recycling rate as well as the buy-back price offered by the manufacturer are higher than those when the government conducts reward-penalty mechanism exclusively for carbon emission; (2) when the government implements carbon emission constraint，both retailers' selling prices of the new product are higher than those when no carbon emission constraint is forced; (3) there is no certain relationship between the two retailers' selling prices of the new product when the government implements the reward-penalty mechanism only for carbon emission and when it implements the mechanism for carbon emission as well as recycling; (4) Regardless of the retailer's fixed cost is high or low, when the government implement reward-penalty mechanism for carbon emission, the retailer's buy-back price is affected by the carbon emission in the manufacturer; (5) When it satisfy certain some conditions, the buy-back price of the retailers' competition reverse supply chain decentralized decision-making model with carbon emission constraint and the government's reward-penalty mechanism can larger than that of the retailers' competition reverse supply chain decentralized decision-making model with carbon emission constraint.

In this paper, we discuss the government's reward and penalty mechanism in the presence of asymmetric information and carbon emission constraint when downstream retailers compete in a reverse supply chain network. However, dual asymmetric information can't be considered in this paper. In the future, we should consider the government's reward and penalty mechanism in the presence of dual asymmetric information and carbon emission constraint when downstream retailers compete in a reverse supply chain network. Moreover, we can discuss the government's reward and penalty mechanism in the presence of asymmetric information and carbon emission constraint when the manufacturers compete in reverse supply chain network.